\documentclass{birkmult}

 \newtheorem{thm}{Theorem}[section]
 \newtheorem{cor}[thm]{Corollary}
 \newtheorem{lem}[thm]{Lemma}
 
 \theoremstyle{definition}
 
 \theoremstyle{remark}

 \numberwithin{equation}{section}

\begin{document}

%
%
%
%
%
%
%
%
%
\title[Groups with normal restriction property]
 {\center Groups with normal restriction property}

  \author[Hung P. Tong-Viet]{\center Hung P. Tong-Viet}

\address{
School of Mathematics\\
University of Birmingham\\
Edgbaston, Birmingham, B15 2TT \\
United Kingdom}

\email{tongviet@maths.bham.ac.uk}

\thanks{This work was completed with the support of University of Birmingham.}



\date{\today}
\dedicatory{}

\begin{abstract}
Let $G$ be a finite group. A subgroup $M$ of $G$ is said to be an \emph{NR-subgroup} if,
whenever $K\unlhd M,$ then $K^G\cap M=K$ where $K^G$ is the normal closure of $K$ in $G.$ Using the Classification of
Finite Simple Groups, we prove that if every maximal subgroup of $G$ is an $NR$-subgroup then $G$ is solvable. This gives a positive answer to a conjecture posed in \cite{ber}.
\end{abstract}

\maketitle

\section{Introduction}

All groups considered are finite.
Let $G$ be a group. Following Berkovich in \cite{ber}, a triple $(G,H,K)$ is said to be {\em
special} in $G$ if $K\unlhd H\leq G$ and $H\cap K^G=K,$ where
$K^G$ is the normal closure of $K$ in $G.$
A subgroup $H$ is called an
$NR-$subgroup (Normal Restriction) if, whenever $K\unlhd H,$ then $(G,H,K)$ is special in $G.$
The main result of this paper is a proof of Conjecture $2$ raised in \cite{ber}.
\begin{thm}\label{main}\emph{(\cite{ber} Conjecture $2$)} If all maximal subgroups of $G$ are $NR$-subgroups then $G$ is solvable.\end{thm}
 In order to prove Theorem \ref{main}, we need a result on the factorization of almost simple groups. Unfortunately, we cannot avoid using the  Classification of Finite Simple Groups in the proof of that result (see Theorem \ref{factorization}). Recall that a group $G$ is said to be \emph{almost simple} if $S\unlhd G\leq Aut(S)$ for some non-abelian simple group $S.$ If $K$ is a proper subgroup of $G$ and $H$ is a subgroup of $G$ with $K\leq H< G,$ then $H$ is called a \emph{ proper over-group} of $K$ in $G.$ Moreover, a subgroup $K$ of $G$ is said to be \emph{$p$-local} in $G$  if $K=N_G(P)$ for some non-trivial
$p$-subgroup $P$ of $G,$ $p$  prime. We also say that  $K$  is \emph{local in}  $G$ if $K$ is $p$-local in $G$ for some prime $p.$ Finally, a subgroup $K$ of $G$
is said to be {\em local maximal} if it is both maximal and local in $G.$
\begin{thm}\label{factorization} Let $S$ be a non-abelian simple group
and $S\unlhd G\leq Aut(S).$ Then there exists a non-trivial subgroup $K$ of $S$ such that
all proper over-groups of $K$ in $S$ are local in $S$ and $G=N_G(K)S.$
\end{thm}
The following corollary is used to show that the minimal
counter-example to Theorem \ref{main} is not simple.
\begin{cor}\label{cor1}  Let $S$ be a non-abelian simple
group. Then $S$ contains a local maximal subgroup.
\end{cor}
\begin{proof} Let $G=Aut(S)$ and $K$ be the subgroup of $S$ obtained from Theorem
\ref{factorization}. Consider the set $\mathcal{A}$ of all proper over-groups of
$K$ in $S.$ Clearly, $\mathcal{ A}$ is non-empty and every element
of $\mathcal{ A}$ is a local subgroup of $S$ containing $K.$ The
maximum element of $\mathcal{ A}$ is a maximal subgroup of
$S$ and is local.
\end{proof}

\smallskip


\section{Preliminaries}
In this section, we collect some  results that we need for the proofs of the theorems above.
\begin{lem}\label{lem1}  Let $K\unlhd H \leq G.$ If $H$ is an
$NR$-subgroup of $G$ then $HK^G/K^G$ is an $NR$-subgroup of
$G/K^G.$ In particular, if $K\unlhd G$ and all maximal subgroups of $G$ are $NR$-subgroups,  then
all maximal subgroups of $G/K$ are also $NR$-subgroups.
\end{lem}
\begin{proof} The first statement is Lemma $4(c)$ in \cite{ber}. The second statement follows easily.\end{proof}

\begin{thm}\emph{(\cite{gor} Theorem $4.3$)}\label{th1} Let $P$ be a $p$-Sylow subgroup of a group $G.$
If $P$ lies in the center of $N_G(P)$ then $G$ has a normal $p$-complement.
\end{thm}

\begin{thm}\label{th2}\emph{(\cite{ber} Proposition $7$)}  Let $H$
be a maximal solvable subgroup of $G.$ If $H$ is an $NR$-subgroup
of $G$ then $H=G.$
\end{thm}


\section{Proofs of the Theorems}

{\it Proof of Theorem \ref{factorization}.}
Without loss of generality, we can assume that $G=Aut(S).$
By the Classification of Finite
Simple Groups, if $S$ is a non-abelian simple group then $S$ is a
finite simple group of Lie type, an alternating group of degree at
least $5$ or one of $26$ sporadic groups. In this proof, we treat
the Tits group, ${}^2F_4(2)'$ as sporadic group rather than a group
of Lie type, and in view of the isomorphisms $A_6\simeq L_2(9),\mbox{ and }
A_5\simeq L_2(5),$ we consider $A_5, A_6$ to be groups of Lie
type.

(i) $S$ is a finite simple group of Lie type in characteristic
$p$, $S\not= {}^2F_4(2)'.$
By Proposition $8.2.1$ and Theorem $13.5.4$ in \cite{car}, $S$ has a (B,N)-pair.
Let $B$ be a Borel subgroup of $S.$ Then $B=N_S(U),$  where $U$ is a $p$-Sylow subgroup of $S.$ For any $\theta\in G,$ as $S\unlhd G,$
$U^\theta\leq S^\theta=S,$ and hence $U^\theta$ is a $p$-Sylow subgroup of $S.$ By Sylow's Theorem
$U^\theta=U^g$ for some $g\in S.$
Observe that $$B^\theta=N_S(U^\theta)=N_S(U^g)=B^g.$$  Thus $\theta g^{-1}\in N_G(B),$ so that $\theta\in N_G(B)S,$
and hence $G=N_G(B)S.$
Moreover, if $H$ is any proper
over-group of $B$ in $S,$ then $H$ is a parabolic subgroup of $S$ and $H<S,$
so that $H$ is  $p$-local  in $S.$ Therefore we can choose $K$ to be a Borel subgroup of $S.$

(ii) $S$ is an alternating group of degree $n\geq 7.$ In this case $G=S_n.$ Let $H=S_{n-3}\times
S_3$ and $K=H\cap S.$ Since $n-3>3,$ it follows from
\cite{lps} that $K$ is a maximal subgroup of $S,$ $H$ is a maximal
subgroup of $G,$ and hence $G=HS.$ As $[H:K]=2,$ we have $H=N_G(K),$ so that
$G=N_G(K)S.$ The subgroup $K$ satisfies the
Theorem since it is $3$-local and maximal in $S.$

(iii) $S$ is sporadic or $S= {}^2F_4(2)'.$\\ By \cite{atlas},
$[G:S]=1$ or $2.$ If $G=S$ then we can choose $K$ to be any local
maximal subgroup of $S.$ The pairs $(S,K)$ are given in Table
\ref{Ta1}. Otherwise, as in $(ii),$ choose $H$ to be a maximal
subgroup of $G$ such that $K=H\cap S$ is a local maximal subgroup
of $S.$ Then $K$ will satisfy the conclusion of the
Theorem. The triple $(S,K,H)$ are given in Table ~\ref{Ta2}. The proof is now completed.
$\Box$\\

\begin{table}

 \begin{center}
 \caption{$|Out(S|)=1$}\label{Ta1}
  \begin{tabular}{|c|c|c|c|c|c|c|c|}
   \hline
   S  & $M_{11}$     &$J_1$&$M_{23}$&$M_{24}$&$Ru$&$Co_3$&$Co_2$\\ \hline
   $K$&$2\:\dot{}S_4$&$7:6$&$23:11$&$2^4:A_8$&$5:4\times A_5$&$2\times M_{12}$&$2^{10}:M_{22}:2$
   \\ \hline
   S &Ly        &Th       & $Fi_{23}$       &$Co_1$         &$J_4$  &B&M \\ \hline
   K&$37:18$& $31:15$ &$2\:\dot{} Fi_{22}$&$S_3\times A_9$&$37:12$&$47:23$&$2\:\dot{} B$\\ \hline
  \end{tabular}

   \end{center}
\end{table}

\begin{table}
 \begin{center}
  \caption{$|Out(S)|=2$}\label{Ta2}
  \begin{tabular}{|c|c|c|}
   \hline
   S  & K & H\\ \hline
   $M_{12}$&$4^2:D_{12}$&$K\cdot 2$\\ \hline
   $M_{22}$&$2^4:A_6$&$2^4:S_6$\\ \hline
   $J_2$&$A_4\times A_5$&$K:2$\\ \hline
   ${}^2F_4(2)'$&$5^2:4A_4$&$5^2:4S_4$\\ \hline
   $HS$&$5:4\times A_5$&$5:4\times S_5$\\ \hline
   $J_3$&$2_{-}^{1+4}:A_5$&$2_{-}^{1+4}:S_5$\\ \hline
   $McL$&$5_+^{1+2}:3:8$&$K\cdot 2$\\ \hline
   $He$&$5^2:4A_4$&$5^2:4S_4$\\ \hline
   $Suz$&$3^5:M_{11}$&$3^5:(M_{11}\times 2)$\\ \hline
   $O'N$&${4^3}^.L_3(2)$&${4^3}^.(L_3(2)\times 2)$\\ \hline
   $Fi_{22}$&$2^{10}:M_{22}$&$2^{10}:M_{22}:2$\\ \hline
   $HN$&$3_+^{1+4}:4A_5$&$3_+^{1+4}:4S_5$\\ \hline
   $Fi_{24}'$&${3^7}^. O_7(3)$&${3^7}^. O_7(3):2$\\ \hline
  \end{tabular}

 \end{center}
\end{table}

\noindent
{\it Proof of Theorem \ref{main}.}
Let $G$ be a minimal
counter-example to Theorem \ref{main}. We first show that $G$ is not simple. By
contradiction, suppose that $G$ is simple. By Corollary
\ref{cor1}, $G$ contains a $p$-local maximal subgroup $M.$ Let
$P$ be a $p$-subgroup of $G$ such that $M=N_G(P).$ Then $1\not=P\unlhd M$ and since
$M$ is an $NR$-subgroup of $G,$ we have $P^G\cap M=P.$ However as  $G$ is simple
and $P\leq P^G\unlhd G,$ $P^G=G.$ Hence $P=G\cap M=M.$ Let $P_1$
be a cyclic subgroup of order $p$ in the center of  $M.$ Then $P_1$ is normal in $M.$
Apply the same argument  as above, we have
$P_1^G=G,$ and so $P_1=P_1^G\cap M=M.$ Thus $M$ is
a cyclic group of order $p.$ In view of the maximality of $M$ and the simplicity of $G,$ $M$ is a $p$-Sylow subgroup
 of $G$ and $N_G(M)=M.$
By Theorem \ref{th1}, $G$ has a normal $p$-complement. This contradicts to our assumption.
Thus $G$ is not simple.

Let $N$ be any minimal normal subgroup of $G.$ By
Lemma \ref{lem1}, the group $G/N$ satisfies the hypothesis of the Theorem and has smaller order than
that of $G,$ by the minimality of $G,$ $G/N$ is solvable. Thus $N$
is the unique minimal normal subgroup of $G,$ and it coincides with
the last term of the derived series of $G.$ If $N$ is solvable
then $G$ is also solvable and we are done. Thus we assume that $N$
is not solvable. Then $N=S_1\times S_2\times\cdots \times S_t,$
where $S_i= S^{x_i},S$ is a non-abelian simple group, and
$x_1, x_2,\cdots,x_t\in G.$ Let $K$ be the subgroup of $S$
obtained from Theorem \ref{factorization}, and $T=K_1\times K_2\times\cdots
\times K_t,$ where $K_i=K^{x_i}.$  Then $T$ is a non-trivial
proper subgroup of $N.$ Since $N$ is the unique minimal normal
subgroup of $G,$ $N_G(T)<G.$ We will show that $G=N_G(T)N.$ For
any $g\in G,$ since $N^g=N,$ there exists a permutation $\pi$ of degree $t$ acting on $\{1,2,\cdots, t\}$ such that $S^{x_ig}=S^{x_{i\pi}}.$ Let $g_i=x_i g
x_{i\pi}^{-1}.$ Then $g_i\in N_G(S).$ We have
$$T^g=K^{x_1g}\times K^{x_2g}\times\cdots \times
K^{x_tg}=K^{g_1x_{1\pi}}\times K^{g_2x_{2\pi}}\times\cdots \times
K^{g_tx_{t\pi}}=$$ $$=K^{g_{1\pi^{-1}}x_1}\times
K^{g_{2\pi^{-1}}x_2}\times\cdots \times
K^{g_{t\pi^{-1}}x_t}=K^{h_1x_1}\times K^{h_2x_2}\times\cdots
\times K^{h_{t}x_t}=$$ $$=K^{x_1s_1}\times K^{x_2s_2}\times\cdots
\times K^{x_ts_t}=K_1^{s_1}\times K_2^{s_2}\times\cdots\times K_t
^{s_t},$$ where $K^{g_{i\pi^{-1}}}=K^{h_i}$ with $h_i\in S$ by
Theorem \ref{factorization},  and $s_i=h_i^{x_i}\in S_i.$ Let
$s=s_1.s_2\dots s_t\in N.$ Since $[S_i,S_j]=1$ if $i\not=j
\in\{1,2,\dots, t\},$ $K_i^s=K_i^{s_i}.$ Thus $T^g=T^s, $ where
$s\in N.$ Therefore $G=N_G(T)N.$

Let $M$ be any maximal subgroup of $G$ containing $N_G(T).$ Let $U=M\cap N.$
We have $G=MN,$ and $U=M\cap N\unlhd M.$ As $G/N=MN/N\simeq M/U,$
$M/U$ is solvable. If $U$ is solvable then $M$ is solvable. By
Theorem \ref{th2}, $G=M,$ a contradiction. Thus $U$ is non-solvable.
Let $L$ be any non-trivial normal subgroup of $M.$ Since $M$ is maximal in $G,$
$M$ is an $NR$-subgroup  of $G,$ so that
$L=L^G\cap M.$ It follows from the fact that $N$ is the
unique minimal normal subgroup of $G,$ $N\leq L^G.$ We have $U=N\cap
M\leq L^G\cap M=L.$ We conclude that $U$ is a minimal normal
subgroup of $M.$ Now, since $U$ is a minimal normal subgroup of $M$ and $U$ is non-solvable, $U=W_1\times
W_2\times \cdots\times W_k,$ where $W_i\simeq W$ for all $1\leq i\leq k$ and $W$ is a non-abelian
simple group. Suppose that there exists $j\in \{1,2,\dots, t\}$
such that $S_j\leq U.$ As $S_j$ is normal in $N,$
$$S_j^G=S_j^{NM}=S_j^M\leq M.$$ However  as $S_j^G=N,$
 $G=MN=M,$ a contradiction. Therefore $S_j\cap U<S_j$ for any $j\in
\{1,2,\dots, t\}.$ Since $K_j\leq S_j\unlhd N,$ $K_j\leq S_j\cap
U\unlhd U.$ As $U$ is a direct product of non-abelian simple groups and $S_j\cap U$
is a non-trivial normal subgroup of $U,$ there exists a non-empty set $J\subseteq
\{1,2,\dots, t\}$ such that $S_j\cap U=\prod_{i\in J}W_i.$ Hence
$$K_j\leq \prod_{i\in J}W_i<S_j,$$ and so $$K\leq \prod_{i\in
J}W_i^{x_j^{-1}}<S,$$ where $W_i^{x_j^{-1}}$ are non-abelian simple
for any $i\in J.$ However, by Theorem \ref{factorization}, $\prod_{i\in
J}W_i^{x_j^{-1}}$ is local in $S.$ This final
contradiction completes the proof.    $  \Box$


\subsection*{Acknowledgment}
I would like to thank Professor Kay Magaard and Doctor Le Thien Tung for their help with the preparation of this work. I am also grateful to Professor Derek Holt and the referee for their suggestions to improve the proof of Theorem \ref{factorization}.
\end{document}